\documentstyle{amsart}
\input{bull-art}
\newcommand{\ZZ}{{\Bbb Z}}
\newcommand{\af}{\alpha}
\newcommand{\dd}{\ldots}
\newcommand{\dr}{\Rightarrow}
\newcommand{\eeq}{\end{equation}}
\newcommand{\beql}[1]{\begin{equation}\label{#1}}
\newcommand{\eqn}[1]{(\ref{#1})}
\newcommand{\hsp}{\hspace*{\parindent}}

\theoremstyle{plain}
\newtheorem{theo}{Theorem}
\newtheorem{thm}{Theorem 7}

\theoremstyle{remark}
\newtheorem{rem}{Remark}

\begin{document}
\def\currentvolume{29}
\def\currentissue{2}
\def\currentyear{1993}
\def\currentmonth{October}
\def\copyrightyear{1993}
\def\currentpages{218-222}
\title {A Linear Construction for Certain Kerdock \\
and Preparata Codes}
\author[A. R. Calderbank, et al.]{A.~R. Calderbank}
\address[A. R. Calderbank and N. J. A. Sloane]
{Mathematical Sciences Research Center \\
AT\&T Bell Laboratories, Murray Hill, New Jersey 07974}
\email[A. Calderbank]{rc@@research.att.com\newline\indent
{\it E-mail address},  N. Sloane: njas@@research.att.com}
\author[]{A. R. Hammons, Jr.}
\address[A. R. Hammons, Jr.]
{Hughes Aircraft Company \\
Canoga Park, California 91304}
\email{hammons@@solar.usc.edu}
\author[]{P. Vijay Kumar}
\address[P. Vijay Kumar]
{Communication Science Institute, EE-Systems \\
University of Southern California, Los Angeles, California 
90089}
\email{kumar@@lamarr.usc.edu}
\author[]{N.~J.~A. Sloane}
\author[]{Patrick Sol\'{e}} 
\address[Patrick Sol\'e]
{Centre National de la
Recherche Scientifique, Sophia -- Antipolis, 06560 
Valbonne, France }
\email{sole@@mimosa.unice.fr} 
\date{November 3, 1992 and, in revised form, December 7, 
1992}
\subjclass{Primary 94B05, 94B15, 94B60}

\maketitle
\begin{abstract}
The Nordstrom-Robinson, Kerdock, and (slightly modified) 
Pre\-
parata codes are
shown to be linear over $\ZZ_4$, the integers $\bmod~4$.
The Kerdock and Preparata codes are duals over $\ZZ_4$,
and the Nordstrom-Robinson code is self-dual.
All these codes are just extended cyclic codes over $\ZZ_4$.
This provides a simple definition for these codes
and explains why their Hamming weight distributions are 
dual to each other.
First- and second-order Reed-Muller codes are also linear 
codes over
$\ZZ_4$, but Hamming codes in general are not,
nor is the Golay code.
\end{abstract}

\section{Introduction}

Some of the best-known examples of nonlinear
binary error-correcting codes that are better than any 
linear codes are
the Nordstrom-Robinson, Kerdock, and Preparata codes
\cite{Ke72, NR67, Pr68, MS77}.
Besides their excellent error-correcting capabilities, 
these codes
are remarkable because the Kerdock and Preparata codes are
``formal duals'', in the sense that although these codes 
are nonlinear,
the weight distribution of one is the MacWilliams
transform of the weight distribution of the other 
\cite[Chapter 15]{MS77}.
The main unsolved question concerning these codes has 
always been whether
they are duals in some more algebraic sense.
Many authors have investigated these codes and have found 
that
(except for the Nordstrom-Robinson code)
they are not unique and, indeed, that large numbers of 
codes exist
with the same weight distributions
\cite{BLW83, Ca89, Ka82, Ka82a, Ka83, VL83}.
Kantor \cite{Ka83} declares that the ``apparent 
relationship between these
[families of codes] is merely a coincidence''.

Although this may be true for many versions of these codes,
we will show that, when properly defined,
Kerdock and Preparata codes are {\em linear} over
$\ZZ_4$ (the integers $\bmod~4$) and that as $\ZZ_4$-codes
they {\em are} duals.
All these codes are, in fact, just extended cyclic codes.

The version of the Kerdock code that we use is the 
standard one,
while our version of the Preparata code differs from the 
standard
one in that it is not a subcode of the Hamming code but of a
nonlinear code with the same weight distribution as the 
Hamming code.
Since the new construction is so simple, we propose that 
this is the
``correct'' way to define these codes.

Kerdock and Preparata codes exist for all lengths $n=4^m 
\ge 16$.
At length 16 they coincide,
giving the Nordstrom-Robinson code \cite{NR67}.
The $\ZZ_4$ version of the Nordstrom-Robinson code is the
``octacode'', a self-dual code of length 8 over $\ZZ_4$ 
that is obtained when the Leech lattice is decomposed into 
eight copies of
the face-centered cubic lattice.
This result was announced in \cite{FST93}.
(The octacode itself is described in \cite{CS92, CS93} and
\S3.)

The very good nonlinear codes of minimal distance 8 
discovered by
Goethals \cite{Go74, Go76} and the high minimal distance 
codes of
Delsarte and Goethals \cite{DG75}
also have a simple description as codes over $\ZZ_4$ (see 
\cite{HKCSS}).

This work developed out of the discovery that four-valued 
sequences have
excellent correlation properties \cite{So89, Bo90, BHK92}) 
and
was carried out independently by Hammons and Kumar 
\cite{HK93} and (very
slightly later) by the other three authors. Theorems 4--6
appear in Hammons's dissertation \cite{Ha92}. In view of the
considerable overlap we have now joined forces. This 
announcement is a
compositum of our results, and full details will be given 
in \cite{HKCSS}.

For undefined terminology from coding theory see 
\cite{MS77}.

\section{Codes over $\ZZ_4$}

A quaternary linear code is an additive subgroup of 
$\ZZ_4^n$.
Duality is defined with respect to the inner product
$a \cdot b = a_1 b_1 + \cdots + a_n b_n$ $(\bmod~4)$.
We define three maps from $\ZZ_4$ to $\ZZ_2$ by
$$
\begin{array}{cccc}
i & \af (i) & \beta (i) & \gamma (i) \\
0 & 0 & 0 & 0 \\
1 & 1 & 0 & 1 \\
2 & 0 & 1 & 1 \\
3 & 1 & 1 & 0
\end{array}
$$
Then we construct binary codes from quaternary codes using 
the map
$\phi: \ZZ_4^n \to \ZZ_2^{2n}$ given by
\begin{equation}
\label{eq1}
\phi (a) = ( \beta (a) , \gamma (a)) ~.
\end{equation}
A binary code is $\ZZ_4$-{\em linear} if its
coordinates can be permuted so that it is the image under 
this
map of a linear code over $\ZZ_4$.

\begin{theo}
\label{th1}
\RM{(a)}~The binary image $\phi (D)$ of a quaternary 
linear code $D$ is 
linear if and only if
\begin{equation}
\label{eq2}
a,b \in D \dr 2 \af (a) \ast \af (b) \in D ~,
\end{equation}
where $\ast$ is componentwise multiplication.

\RM{(b)}~A binary linear code $C$ of even length is 
$\ZZ_4$-linear if and only if
its coordinates can be permuted so that
\begin{equation}
\label{eq3}
u, v \in C \dr
(u+s(u)) \ast (v+s(v)) \in C ~,
\end{equation}
where $s$ is the \RM{``}swap\RM{''} 
map that interchanges the left and right halves of a vector.
\end{theo}

\begin{theo}
\label{th2}
Binary Reed-Muller codes of length $n=2^m \ge 2$ and 
orders $0,1,2,
 m-1, m$
are $\ZZ_4$-linear.
\end{theo}

\begin{theo}
\label{th3}
Extended Hamming codes of lengths $n=2^m \ge 32$ are not 
$\ZZ_4$-linear, nor is the Golay code of length 
\RM{24}.
\end{theo}

\begin{rem}
In \eqn{eq3} if $u,v$ are represented by Boolean functions 
of degree $r$
and $(u+s(u)) \ast (v+s(v)) \neq 0$, then
$(u+s(u)) \ast (v+s(v))$ is a Boolean function of degree 
$2r-2$.
So an $r$th-order Reed-Muller code with $r \le m/2$ 
satisfies
\eqn{eq3} provided $r \le 2$ and we conjecture it
does not satisfy \eqn{eq3} if $3 \le r \le m-2$.
The first assertion of Theorem~\ref{th3} establishes that 
$(m-2)$nd-order Reed-Muller codes are not
$\ZZ_4$-linear for $m \ge 5$.
\end{rem}

Let $D$ be a quaternary linear code and $C = \phi (D)$ the 
corresponding
binary code.
In general $C$ is not linear,
but we define the $\ZZ_4$-dual of $C$ to be $C^{\perp_4} = 
\phi (D^\perp )$,
where $D^\perp$ denotes the dual  code to $D$,
as in the following diagram.
$$
\begin{array}{r@{~}ccc@{~}c@{~}c}
~ & D & \stackrel{\phi}{\longrightarrow} & C & = & \phi 
(D) \\
\mbox{dual} &  \downarrow \\
~ &D^\perp & \stackrel{\phi}{\longrightarrow} & 
C^{\perp_4} & = & \phi (D^\perp )
\end{array}
$$

The familiar Hamming weight enumerator for a binary linear 
code $C$
will be denoted by $W_C (x,y)$.
This weight enumerator is also well defined for binary 
nonlinear codes provided
they are distance-invariant \cite{MS77}.
The {\em symmetrized weight enumerator}
of a quaternary linear code $D$ is
$$
\mathrm{swe}_{D} (x,y,z) = \sum_{a \in D} x^{N_0 (a)} 
y^{N_1 (a)} z^{N_2 (a)} ~,$$
where $N_i (a)$ is the number of components of $a$ congruent
to $\pm i$ $(\bmod~4)$.
Then (see \cite{CS93, Kl87})
\[\label{eq4}
\mathrm{swe}_{D^\perp}
(x,y,z) = \frac{1}{|D|}
~\mathrm{swe}_D (x+2y+z , x-z , x-2y +z ) ~.
\]

\begin{theo}
\label{th4}
If $D$ is a quaternary linear code, then $C = \phi (D)$ and
$C^{\perp_4} = \phi (D^\perp )$ are distance invariant, and
\[\label{eq5}
W_C (x,y) = \mathrm{swe}_D (x^2 , xy , y^2) ~,
\]

\[\label{eq6}
W_{C^{\perp_4}} (x,y) = \frac{1}{|C|}
W_C (x+y , x-y) ~.
\]
\end{theo}
\section{Kerdock, Preparata and Nordstrom-Robinson codes}
\hsp
Let $h_2 (X) \in \ZZ_2 [X]$ be a primitive irreducible 
polynomial of degree $m$.
There is a unique monic polynomial
$h(X) \in \ZZ_4 [Z]$ of degree $m$ such that $h(X) \equiv 
h_2 (X)$
$(\bmod~2)$ and $h(X)$ divides $X^n -1$ $(\bmod~4)$,
where $n=2^m -1$
\cite{So89, Ya90}.
Let $g(X)$ be the reciprocal polynomial to 
$$(X^n -1 ) / ((X-1) h(X)).$$

\begin{theo}
\label{th5}
The cyclic code generated by $g(X)$, extended by an overall
parity check, is a quaternary code $D$ of length $2^m$
containing $4^{m+1}$ words.
For $m$ odd $\ge 3$ the corresponding binary code
$K = \phi (D)$ is the Kerdock code of length $2^{m+1}$ 
containing
$2^{2m+2}$
words and with minimal distance $2^m -2^{(m-1)/2}$.
\end{theo}

The proof is by showing that $D$ has a simple definition 
in terms of
the relative trace function
from $\ZZ_4 [ \xi ]$ to $\ZZ_4$, where $\xi$ is a root of 
$h(X)$
(cf. \cite{Bo90, BHK92}),
and in this form it agrees with the definition of the 
Kerdock
code given in
\cite[pp.~457--458]{MS77}.

\begin{theo}
\label{th6}
The cyclic code generated by $h(X)$, extended by an 
overall parity check,
is a quaternary code $D^\perp$ dual to $D$.
For $m$ odd $\ge 3$
the corresponding binary code $P = \phi (D^\perp )$ has
length $2^{m+1}$, contains $2^k$ words,
$k=2^{m+1} - 2m-2$, has minimal
distance \RM6,
and has the same weight enumerator as the Preparata code.
\end{theo}

$P$ is the $\ZZ_4$-dual of $K$, both codes are distance 
invariant, and the
weight distribution of one is the MacWilliams transform of 
the weight
distribution of the other.

For example when $m=5$, we may take 
$$h(X) = \sum_{i=0}^5 h_i X^i,\qquad
g(X) = \sum_{i=0}^{25} g_i X^i,$$
where $h_0 \dd$ and $g_0 \dd$ are
323001 and 11120122010303133013212213.
In this case the linear span of $P$ has minimal distance 2,
which shows that $P$ is strictly different from 
Preparata's original
construction \cite{Pr68}, for which the
linear span is the extended Hamming code.

There is a distance-regular graph \cite{BCN89} defined on 
cosets of our
Preparata code which may be of some combinatorial interest.

In the case $m=3$, both $P$ and $K$ become the 
Nordstrom-Robinson code
\cite{NR67}, and the quaternary code $D= D^\perp$ is the 
octacode. The latter
may be defined as the extended cyclic code generated by 
$h(X) = X^3 + 3X^2 +
2X+3$ or as the unique quaternary self-dual code of length 
8 which has the
property that its binary image has minimal distance 6 
\cite{CS93} or as the
``glue code'' used to construct the Leech lattice from a 
direct sum of eight
copies of the face-centered cubic lattice $A_3$ (note that 
$A_3^\perp / A_3
\cong \ZZ_4$) \cite[Chapter 24]{CS92}.

\begin{thm}[\cite{FST93}]
The Nordstrom-Robinson code is the binary image of the 
octacode.
\end{thm}


\begin{thebibliography}{MMM991}

\bibitem[BLW83]{BLW83}
R.~D. Baker, J.~H. van~Lint, and R.~M. Wilson,
{\em On the Preparata and Goethals codes},
IEEE Trans. Inform. Theory {\bf 29} (1983), 342--345.

\bibitem[Bo90]{Bo90}
S. Bozta\c{s},
{\em Near-optimal $4 \phi$ \rom{(4}-phase\rom{)} sequences 
and optimal binary sequences for
CDMA},
Ph.D. dissertation, Univ. of Southern California, Los 
Angeles, 1990.

\bibitem[BHK92]{BHK92}
S. Bozta\c{s},
A.~R. Hammons, Jr., and P.~V. Kumar,
{\em $4$-phase sequences with near-optimum correlation 
properties},
IEEE Trans. Inform. Theory
{\bf 38} (1992), 1101--1113.

\bibitem[BCN89]{BCN89}
A.~E. Brouwer,
A.~M. Cohen, and A. Neumaier,
{\em Distance-regular graphs},
Springer-Verlag, New York, 1989.

\bibitem[Ca89]{Ca89}
C. Carlet,
{\em A simple description of Kerdock codes},
Lecture Notes in Comput. Sci., vol. 388, Springer-Verlag,
Berlin and New York, 1989, pp. 202--208.

\bibitem[CS92]{CS92}
J.~H. Conway and N.~J.~A. Sloane,
{\em Sphere-packings, lattices and groups},
2nd ed., Springer-Verlag, New York, 1992.

\bibitem[CS93]{CS93}
\bysame, 
{\em Self-dual codes over the integers modulo} 4,
J. Combin. Theory Ser.~A {\bf 62} (1993), 30--45.

\bibitem[DG75]{DG75}
P. Delsarte and J.~M. Goethals,
{\em Alternating bilinear forms over $GF (q)$},
J. Combin. Theory Ser. A {\bf 19} (1975), 26--50.

\bibitem[FST93]{FST93}
G.~D. Forney, Jr.,
N.~J.~A. Sloane, and M.~D. Trott,
{\em The Nordstrom-Robinson code is the binary image of 
the octacode},
Proceedings DIMACS/IEEE Workshop on Coding and Quantization,
DIMACS Series in Discrete Mathematics and
Theoretical Computer Science, Amer. Math. Soc., 
Providence, RI
(to appear).

\bibitem[Go74]{Go74}
J.~M. Goethals,
{\em Two dual families of nonlinear binary codes},
Electron. Lett. {\bf 10} (1974), 471--472.

\bibitem[Go76]{Go76}
\bysame,
{\em Nonlinear codes defined by quadratic forms over 
$GF(2)$},
Inform. Control {\bf 31} (1976), 43--74.

\bibitem[Ha92]{Ha92}
A.~R. Hammons, Jr.,
{\em On four-phase sequences with low correlation and 
their relation
to Kerdock and Preparata codes},
Ph.D. dissertation, Univ. of Southern California,
November 1992.

\bibitem[HK93]{HK93}
A.~R. Hammons, Jr., and P.~V. Kumar,
{\em On the apparent duality of Kerdock and Preparata 
codes},
Abstracts, IEEE Internat. Sympos. Inform. Theory,
San Antonio, TX, January 1993.

\bibitem[HKCSS]{HKCSS}
A.~R. Hammons, Jr., P. V. Kumar, A.~R. Calderbank, 
N.~J.~A. Sloane, and
P. Sol\'{e},
{\em The $\ZZ_4$-linearity of Kerdock, Preparata, Goethals 
and related
codes},
IEEE Trans. Inform. Theory, in press.

\bibitem[Ka82]{Ka82}
W.~M. Kantor,
{\em An exponential number of generalized Kerdock codes},
Inform. Control {\bf 53} (1982), 74--80.

\bibitem[Ka82a]{Ka82a}
\bysame,
{\em Spreads, translation planes and Kerdock sets},
SIAM J. Algebra Discrete Math. {\bf 3} (1982),
151--165, 308--318.

\bibitem[Ka83]{Ka83}
\bysame,
{\em On the inequivalence of generalized Preparata codes},
IEEE Trans. Inform. Theory {\bf 29} (1983), 345--348.

\bibitem[Ke72]{Ke72}
A.~M. Kerdock,
{\em A class of low-rate nonlinear binary codes},
Inform. Control {\bf 20} (1972), 182--187.

\bibitem[Kl87]{Kl87}
M. Klemm,
{\em \"Uber die Identit\"{a}t von MacWilliams f\"{u}r die
Gewichtsfunktion von Codes},
Arch. Math. (Brno) {\bf 49} (1987), 400--406.

\bibitem[VL83]{VL83}
J.~H. van~Lint,
{\em Kerdock and Preparata codes},
Congr. Numer. {\bf 39} (1983), 25--41.

\bibitem[MS77]{MS77}
F.~J. MacWilliams and N.~J.~A. Sloane,
{\em The theory of error-correcting codes},
North-Holland, Amsterdam, 1977.

\bibitem[NR67]{NR67}
A.~W. Nordstrom and J.~P. Robinson,
{\em An optimum nonlinear code},
Inform. Control {\bf 11} (1967), 613--616.

\bibitem[Pr68]{Pr68}
F.~P. Preparata,
{\em A class of optimum nonlinear double-error correcting 
codes},
Inform. Control {\bf 13} (1968), 378--400.

\bibitem[So89]{So89}
P. Sol\'{e},
{\em A quaternary cyclic code, and a family of quadriphase
sequences with low correlation properties},
Lecture Notes in Comput. Sci., vol. 388, Springer-Verlag,
New York and Berlin, 1989, pp. 193--201.

\bibitem[Ya90]{Ya90}
M. Yamada,
{\em Distance-regular digraphs of girth \RM4 over an 
extension ring of
$Z/4Z$},
Graphs Combin. {\bf 6} (1990), 381--394.
\end{thebibliography}
\end{document}